\documentclass{article}
\usepackage{amsfonts,amssymb}
\usepackage{stmaryrd}  

\usepackage{color} 
   \definecolor{cites}{rgb}{0.50 , 0.00 , 0.00}  
   \definecolor{urls} {rgb}{0.00 , 0.00 , 0.50}  
   \definecolor{links}{rgb}{0.00 , 0.00 , 0.50}   
\usepackage[
      colorlinks=true,   
      citecolor=cites,   
      urlcolor=urls,     
      linkcolor=links,   
      pdftitle={Approximating the inverse of banded matrices by banded matrices with applications to probability and statistics},  
      pdfauthor={Peter J. Bickel, Marko Lindner}, 
      pdfpagemode=UseOutlines, 
      pdfstartview=FitH,       
      bookmarksopen=false      
   ]{hyperref}

\newcommand\C{{\mathbb C}}
\newcommand\R{{\mathbb R}}
\newcommand\Z{{\mathbb Z}}
\newcommand\N{{\mathbb N}}
\newcommand\T{{\mathbb T}}
\newcommand\I{{\mathcal I}}
\newcommand\W{{\mathcal W}}
\newcommand\Ta{{\mathcal T}}

\newcommand\BL{{\rm BL}}
\newcommand\BO{{\rm BO}}
\newcommand\BDO{{\rm BDO}}
\newcommand\HBO{{\rm HBO}}
\newcommand\HBDO{{\rm HBDO}}

\newcommand\dist{{\rm dist}}
\newcommand\spec{{\rm sp\,}}
\newcommand\toto{{\rightrightarrows}}
\newcommand\lW{{\llbracket}}
\newcommand\rW{{\rrbracket}}

\newcommand\rh{\varrho}

\newtheorem{theorem}{Theorem}[section]
\newtheorem{lemma}[theorem]{Lemma}

\newtheorem{proposition}[theorem]{Proposition}

\newenvironment{remark}
 {\par\noindent\refstepcounter{theorem}{\bf Remark \thetheorem}}
 {\raisebox{1mm}{\framebox{}}\pagebreak[2]}

\newenvironment{example}
 {\par\noindent\refstepcounter{theorem}{\bf Example \thetheorem}}
 {\raisebox{1mm}{\framebox{}}\pagebreak[2]}

\newenvironment{proof}
 {\par\noindent{\it Proof.}}
 {\rule{2mm}{2mm}\pagebreak[2]}

\newcommand{\be}{\begin{equation}}
\newcommand{\ee}{\end{equation}}
\newcommand{\ve}{\varepsilon}
\newcommand{\X}{\mathbf{X}}
\newcommand{\Y}{\mathbf{Y}}

 \newcommand{\mc}{\mathcal}
\newcommand{\0}{\mbox{\boldmath$0$}}
\newcommand{\bsig}{\mbox{\boldmath$\sigma$}}

\newcommand{\bx}{{\bf x}}

\newcommand{\bt}{{\bf t}}
\newcommand{\bn}{{\bf n}}

\newcommand\eqref[1]{(\ref{#1})}

\begin{document}
\title{\bf Approximating the inverse of banded matrices\\
by banded matrices with applications\\
to probability and statistics}
\author{{\sc Peter Bickel}\ \ and\ \ {\sc Marko Lindner}}
\date{\today}
\maketitle
\begin{quote}
\renewcommand{\baselinestretch}{1.0}
\footnotesize {\sc Abstract.} In the first part of this paper
we give an elementary proof of the fact that if an infinite
matrix $A$, which is invertible as a bounded operator on $\ell^2$,
can be uniformly approximated by banded matrices then so can the
inverse of $A$. We give explicit formulas for the banded
approximations of $A^{-1}$ as well as bounds on their accuracy
and speed of convergence in terms of their band-width.
In the second part we apply these results to covariance matrices
$\Sigma$ of Gaussian processes and study mixing and beta mixing of
processes in terms of properties of $\Sigma$. Finally, we note some
applications of our results to statistics.
\end{quote}


\section{Introduction}
Let $\I$ be either the set of natural numbers, $\N$, or the
integers, $\Z$, and let $\ell^2$ denote the corresponding Hilbert
space of one- or two-sided infinite sequences $(u_k)_{k\in\I}$ of
complex numbers with $\sum_{k\in\I}|u_k|^2<\infty$. We know that
every bounded linear operator $A$ on $\ell^2$ can naturally be
identified with a (one- or two-sided) infinite matrix
$(a_{ij})_{i,j\in\I}$. We will therefore use the words 'operator'
and 'matrix' synonymously here.

It is clear that the inverse of a banded infinite matrix $A$, if it
exists, is in general not a banded matrix any more. However, one can
show that it is still the uniform limit of a sequence of such, that
is: it is what we call a {\sl band-dominated matrix}. By a simple
approximation argument, this result immediately implies that the
class of band-dominated matrices is {\sl inverse closed}, that
means: if one of them is invertible, its inverse is again
band-dominated.

{\bf Outline of the Paper. } We will give a proof of this result in
Section \ref{sec:BDO} which moreover comes with explicit formulas
for the banded approximations of the inverse and with bounds on
their accuracy and speed of convergence in terms of their
band-width. In Section \ref{sec:Wiener} we apply these formulas to
another class of operators, the so-called {\sl Wiener algebra}.
These inverse closedness results themselves are not new (see e.g.
Kurbatov \cite{Kurbatov82,Kurbatov89,KurbatovBook} but also
\cite{LaRa82,LevRa79,LutRa78,Shubin78,Shubin86,RaRoSi:Book,Li:Book}
for related questions) but what we believe is new here is our
approach and the explicit approximates of the inverse that it comes
with, as well as the generalizations of our results to the operator
classes defined in Section \ref{sec:gen} with an eye to applications
in statistics.

Next, in Section \ref{sec:Gauss}, we study the
relation between our results in Sections \ref{sec:BDO} and
\ref{sec:Wiener} and the notion of {\sl regularity} in Gaussian
processes, a well settled problem in the stationary case, clarified,
as we believe, for the first time in the general (non-stationary)
case. In the same section, we consider the characterization of the
notion of {\sl beta mixing} (or {\sl absolute regularity}) for Gaussian
processes, a problem considered by Ibragimov and Solev \cite{IbragSolev}
for the stationary case. We give a necessary and sufficient condition for
beta mixing in the general case, relating it to closure notions
stronger than those in the previous sections. Finally, in Section
\ref{sec:stat}, we sketch the applications to statistics which
initially prompted this work.

~

{\bf Notations. }
Here and in what follows, let $\ell^2:=\ell^2(\I)$ stand for
the set of all sequences $x=(x_k)_{k\in\I}$ of complex numbers with
\[
\|x\|\ :=\ \sqrt{\sum_{k\in\I} |x_k|^2}\ <\ \infty,
\]
where the index set $\I$ is fixed. For every bounded and linear
operator $A$ on $\ell^2$, let $A^*$ denote the adjoint operator with
matrix representation $(a_{ij})^*=(\overline{a_{ji}})$ and let
$\|A\|$ denote the induced operator norm, that is
\[
\|A\|\ =\ \sup_{x\in\ell^2\setminus\{0\}}\frac{\|Ax\|}{\|x\|}.
\]
By putting $\I=\{1,...,n\}$ in the above, we define $\|x\|$ and
$\|A\|$ as is usual for finite vectors $x\in\C^n$ and $n\times n$
matrices $A$ acting on them. We will now suppose that
$\I\in\{\Z,\N\}$.

Let $\BL:=\BL(\ell^2)$ be the set of all bounded linear operators
$A:\ell^2\to\ell^2$. Equipped with addition, multiplication by
scalars, operator composition and the above norm, $\BL$ is a Banach
algebra and even a $C^*$-algebra with involution $A\mapsto A^*$. If
$A\in\BL$ is invertible (i.e. a bijection $\ell^2\to\ell^2$) then
also $A^{-1}\in\BL$ as a consequence of the open mapping theorem.
Now let $\BO:=\BO(\ell^2)$ refer to the set of all operators $A\in
\BL$ that are induced by a banded matrix -- meaning a matrix with
only finitely many nonzero diagonals. Clearly, the set $\BO$ is
closed under addition, multiplication, multiplication by scalars and
under passing to the adjoint -- but it is not closed in the operator
norm $\|\cdot\|$ on $\BL$. That is why one is interested in the norm
closure of $\BO$, henceforth denoted by $\BDO:=\BDO(\ell^2)$, the
elements of which are called {\sl band-dominated
operators/matrices}.

From a computational point of view the operator norm $\|\cdot\|$ is
not very handy. An alternative norm $\lW\cdot\rW$ can be defined on
$\BO$ as follows: For $A\in\BO$ with matrix representation
$(a_{ij})_{i,j\in\I}$ and for each $k\in\Z$, let $d_k$ be the
supremum norm of the $k$-th diagonal of $A$, that is
\begin{equation} \label{eq:Wdef}
d_k\ :=\ \sup\{\,|a_{ij}|\ :\ i,j\in\I,\ i-j=k\,\},\qquad
\textrm{and put}\qquad \lW A\rW\ :=\ \sum_{k\in\Z}d_k.
\end{equation}
It is easy to see that this defines a norm on $\BO$ with
$\|A\|\le\lW A\rW$ for all $A\in\BO$. Let us this time pass to the
completion of $\BO$ in the stronger norm $\lW\cdot\rW$; what we get
is a proper subset of $\BDO$ that shall be denoted by $\W$.
Equivalently, $A\in\W$ iff $\lW A\rW<\infty$, where $\lW A\rW$ is as
defined in (\ref{eq:Wdef}) but now for arbitrary infinite matrices.
It turns out that $(\W,\lW\cdot\rW)$ is a Banach algebra that is
often referred to as the {\sl Wiener algebra}. If we, for a moment,
generalize our setting from operators on $\ell^2$ to operators on
$\ell^p$ with $p\in[1,\infty]$, it is clear that the class $\BO$
does not depend on $p$. Unlike the norm closure $\BDO(\ell^p)$ of
$\BO$, the Wiener algebra $\W$ is also independent of $p$ since it
is defined merely in terms of matrix entries. One has, for all
$p\in[1,\infty]$, that
\begin{equation} \label{eq:WinBDO}
\W\subset\BDO(\ell^p)\qquad\textrm{with}\qquad\|A\|\le\lW A\rW
\end{equation}
for all $A\in\W$. We will give an elementary proof of the inverse
closedness of $\W$, thereby automatically proving Wiener's famous
theorem on functions with absolutely summable Fourier series.

\section{$\BDO$ is inverse closed} \label{sec:BDO}
The shortest proof of the inverse closedness of $\BDO$ goes like
this: By its definition, $\BDO$ is a Banach subalgebra of $\BL$ that
is closed under the involution map $A\mapsto A^*$. Since, moreover,
the equality $\|A^*A\|=\|A\|^2$ holds for all $A\in\BL$, both $\BL$
and $\BDO$ are C$^*$-algebras, and a basic theorem \cite{Dixmier} on
C$^*$-algebras says that therefore $\BDO$ is inverse closed in
$\BL$, i.e. if $A\in\BDO$ is invertible in $\BL$ one always has
$A^{-1}\in\BDO$. In this section we will find out how to approximate
$A^{-1}$ by band matrices -- and how good this approximation is.

In order to distinguish between banded matrices of different {\sl
band widths} $k\in\N_0:=\N\cup\{0\}$, we will introduce the notation
$\BO_k$ for the set of all $A\in\BL$ whose matrix
$(a_{ij})_{i,j\in\I}$ is supported on the diagonals numbered
$-k,...,k$ only, that means $a_{ij}=0$ if $|i-j|>k$, i.e. $d_n=0$ if
$|n|>k$ with $d_n$ from (\ref{eq:Wdef}). Clearly, we have
$\BO_k\subset\BO_{k+1}$, $\BO=\cup_{k\ge 0} \BO_k$ and
\begin{equation} \label{eq:BDOconv}
A\in\BDO\qquad \textrm{iff}\qquad 0\ =\ \dist(A,\BO)\ =\
\dist(A,\bigcup_{k=0}^\infty \BO_k)\ =\ \lim_{k\to\infty}
\dist(A,\BO_k)
\end{equation}
with the usual definition of the distance, $\dist(A,S):=\inf_{B\in
S}\|A-B\|$, of an operator $A\in\BL$ from a set $S\subset\BL$. Note
that if $a_{ij}$ is a matrix entry of $A$ with $|i-j|>k$ then
clearly $a_{ij}$ is still a matrix entry of $A-B$ for all
$B\in\BO_k$ so that $\|A-B\|\ge |a_{ij}|$. Consequently,
\[
\dist(A,\BO_k)\ =\ \inf_{B\in\BO_k}\|A-B\|\ \ge\ |a_{ij}|,\qquad
|i-j|>k
\]
holds, i.e. $\dist(A,\BO_k)$ is a bound on all matrix entries
outside the $-k,...,k$ band of $A$. Using the diagonal suprema
introduced in (\ref{eq:Wdef}), we can rephrase this as
\begin{equation} \label{eq:dist_d}
\dist(A,\BO_k)\ \ge\ \sup\{\,|a_{ij}|\ :\ i,j\in\I,\ i-j=n\,\}\ =\
d_n,\qquad |n|>k.
\end{equation}
We start with the simple case when $A$ is banded and self-adjoint
positive definite, i.e. $A\in\BO$, $A=A^*$ and the spectrum of $A$,
$\spec A$, is strictly positive. In this case, it is well known that
\[
M\ :=\ \sup_{\lambda\in\spec A}|\lambda|\ =\ \varrho(A)\ =\
\|A\|,\qquad m\ :=\ \inf_{\lambda\in\spec A}|\lambda|\ =\
1/\varrho(A^{-1})\ =\ 1/\|A^{-1}\|
\]
with $\varrho(A)$ denoting the spectral radius of $A$. Moreover, we
have that $\kappa:=M/m=\|A\|\|A^{-1}\|$ is the condition number of
$A$.

\begin{lemma} \label{lem:BOpos}
Let $A\in\BO_k$ for some $k\in\N_0$ be self-adjoint positive definite,
and define $M,m$ and $\kappa$ as above. Then, for every $n\in\N_0$,
it holds that
\begin{equation} \label{eq:BOpos}
\dist\big(\,A^{-1}\,,\,\BO_{n\cdot k}\,\big)\ \le\ \frac 1m
\left(\frac{M-m}{M+m}\right)^{n+1}\ =\ \frac 1m
\left(\frac{\kappa-1}{\kappa+1}\right)^{n+1}
\end{equation}
where an approximation $B_n\in\BO_{n\cdot k}$ of $A^{-1}$ with this
accuracy is given in (\ref{eq:Bn}) below. In particular,
$A^{-1}\in\BDO$ since the right-hand side of (\ref{eq:BOpos}) goes
to zero as $n\to\infty$.
\end{lemma}
\begin{proof}
We start by looking for a $\gamma\in\R$ such that
\[
\|I-\gamma A\|\ =\ \varrho(I-\gamma A)\ =\ \sup_{\lambda\in\spec
A}|\gamma\lambda-1|\ =\ \max\big(\,|\gamma m-1|\,,\,|\gamma
M-1|\,\big)
\]
is minimized. A little thought shows that this is the case iff
$|\gamma m-1|=|\gamma M-1|$, i.e. $1$ is the midpoint of the
interval $[\gamma m,\gamma M]$ so that $\gamma=\frac 2{M+m}$. In
this case,
\[
\|I-\gamma A\|\ =\ 1-\gamma m\ =\ 1-\frac {2m}{M+m}\ =\
\frac{M-m}{M+m}\ =\ \frac{\kappa-1}{\kappa+1}\ =\ 1-\frac
2{\kappa+1}\ <\ 1.
\]
Now, by Neumann series, for every $n\in\N_0$,
\begin{equation} \label{eq:Bn}
A^{-1}\ =\ \gamma(\gamma A)^{-1}\ =\ \gamma
\sum_{j=0}^\infty(I-\gamma A)^j\ =\ \underbrace{\gamma
\sum_{j=0}^n(I-\gamma A)^j}_{=:\,B_n}\ +\ \underbrace{\gamma
\sum_{j=n+1}^\infty(I-\gamma A)^j}_{=:\,C_n}
\end{equation}
holds with $B_n\in\BO_{nk}$ and
\begin{equation} \label{eq:Cn}
\|C_n\|\ \le\ |\gamma|\sum_{j=n+1}^\infty\|I-\gamma A\|^j\ =\ \frac
2{M+m}\,\frac{\left(\frac{M-m}{M+m}\right)^{n+1}}{1-\frac{M-m}{M+m}}\
=\ \frac 1m \left(\frac{M-m}{M+m}\right)^{n+1}\ =\ \frac 1m
\left(\frac{\kappa-1}{\kappa+1}\right)^{n+1},
\end{equation}
which finishes the proof.
\end{proof}
\bigskip

We now pass to the non-self-adjoint case -- but still banded.

\begin{proposition} \label{prop:BO}
Let $A\in\BO_k$ for some $k\in\N_0$ be invertible and again put
$M:=\|A\|$, $m:=1/\|A^{-1}\|$ and $\kappa:=M/m=\|A\|\|A^{-1}\|$.
Then, for every $n\in\N_0$, it holds that
\begin{equation} \label{eq:BO}
\dist\big(\,A^{-1}\,,\,\BO_{3nk}\,\big)\ \le\
\frac{M}{m^2}\left(\frac{M^2-m^2}{M^2+m^2}\right)^{n+1}\ =\
\frac{\kappa^2}M\left(\frac{\kappa^2-1}{\kappa^2+1}\right)^{n+1}
\end{equation}
where an approximation of $A^{-1}$ in $\BO_{3nk}$ with this accuracy
is given in (\ref{eq:Bn2}) below. In particular, $A^{-1}\in\BDO$
since the right-hand side of (\ref{eq:BO}) goes to zero as
$n\to\infty$.
\end{proposition}
\begin{proof}
The idea is to write $A^{-1}\ =\ B^{-1}A^*$, where
$B:=A^*A\in\BO_{2k}$ is clearly self-adjoint positive definite, and to
approximate $B^{-1}$ as in the previous lemma. When we apply Lemma
\ref{lem:BOpos} to $B$ (in place of $A$), note that
\begin{eqnarray*}
M_B&:=&\|B\|\ =\ \|A^*A\|\ =\ \|A\|^2\ =\ M^2,\\
m_B&:=&1/\|B^{-1}\|\ =\ 1/\|A^{-1}(A^*)^{-1}\|\ =\
1/\|A^{-1}(A^{-1})^*\|\ =\ 1/\|A^{-1}\|^2\ =\ m^2\quad\textrm{and}\\
\kappa_B&:=&M_B/m_B\ =\ M^2/m^2\ =\ \kappa^2.
\end{eqnarray*}
Now for every $n\in\N_0$, in analogy to (\ref{eq:Bn}), we can write
$B^{-1}=B_n+C_n$ with $B_n\in\BO_{2nk}$ and $\|C_n\|$ bounded as in
(\ref{eq:Cn}), so that $A^{-1}=B^{-1}A^*=B_nA^*+C_nA^*$ with
\begin{equation} \label{eq:Bn2}
B_nA^*\ =\ \gamma_B\sum_{j=0}^n(I-\gamma_BB)^jA^*\ =\ \frac
2{M^2+m^2}\sum_{j=0}^n(I-\frac{2\,A^*A}{M^2+m^2})^jA^*\ \in\
\BO_{3nk}
\end{equation}
and
\[
\|C_nA^*\|\ \le\ \frac 1{m_B}
\left(\frac{M_B-m_B}{M_B+m_B}\right)^{n+1}M\ =
\frac{M}{m^2}\left(\frac{M^2-m^2}{M^2+m^2}\right)^{n+1}\ =\
\frac{\kappa^2}M\left(\frac{\kappa^2-1}{\kappa^2+1}\right)^{n+1},
\]
which proves the result.
\end{proof}
\bigskip

Finally, we pass to the most general case, $A\in\BDO$.

\begin{theorem} \label{th:BDO}
Let $A\in\BDO$ be invertible, put $\delta_k:=\dist(A,\BO_k)$ for
$k=0,1,2,...$ and again let $M:=\|A\|$, $m:=1/\|A^{-1}\|$ and
$\kappa:=M/m=\|A\|\|A^{-1}\|$ as well as
$\alpha_k:=m/(m-2\delta_k)$. Note that, by (\ref{eq:BDOconv}),
$\delta_k\to 0$ and $\alpha_k\to 1$ as $k\to\infty$. Then, for all
$n\in\N_0$ and all sufficiently large $k\in\N_0$ (such that
$\delta_k<m/2$), it holds that
\begin{eqnarray*}
\dist\big(\,A^{-1}\,,\,\BO_{3nk}\,\big)
&\le& \frac{\alpha_k}{m^2}\left(2\delta_k\ +\ \alpha_k(M+2\delta_k)\left(\frac{(M+2\delta_k)^2-(m/\alpha_k)^2}{(M-2\delta_k)^2+(m/\alpha_k)^2}\right)^{n+1} \right)\\
&=& \frac{2\delta_k\alpha_k}{m^2}\ +\
\frac{(\kappa_k^+)^2}{M+2\delta_k}
\left(\frac{(\kappa_k^+)^2-1}{(\kappa_k^-)^2+1}\right)^{n+1}
\end{eqnarray*}
with $\kappa_k^\pm$ defined by (\ref{eq:kappapm}), respectively. In
particular, $A^{-1}\in\BDO$.
\end{theorem}
\begin{proof}
For every $k\in\N_0$ pick an $A_k\in\BO_k$ with $\|A-A_k\|\le
2\delta_k$. Then $A_k\toto A$ (norm convergence in $\BL$) as
$k\to\infty$ since $\delta_k\to 0$. Since $A$ is invertible we know
that $A_k$ is invertible for sufficiently large $k$; precisely, take
$k_0\in\N_0$ big enough that $\delta_k< m/2$ for all $k>k_0$ so that
$\|A^{-1}(A-A_k)\|\le 2\delta_k/m<1$. Then $A_k=A(I-A^{-1}(A-A_k))$
is invertible with
\[
\|A_k^{-1}\|\ =\ \|\sum_{j=0}^\infty(A^{-1}(A-A_k))^jA^{-1}\|\ \le\
\sum_{j=0}^\infty \left(\frac{2\delta_k}m\right)^j \|A^{-1}\|\ =\
\frac 1{1-\frac{2\delta_k}m} \|A^{-1}\|\ =\ \alpha_k\|A^{-1}\|
\]
and hence
\begin{equation} \label{eq:horiz}
\|A_k^{-1}-A^{-1}\|\ =\ \|A^{-1}(A-A_k)A_k^{-1}\|\ \le\
\|A^{-1}\|\,\|A-A_k\|\,\|A_k^{-1}\|\ \le\ \frac 1m\, 2\delta_k\,
\frac{\alpha_k}m\ =\ \frac{2\delta_k\alpha_k}{m^2}
\end{equation}
for all $k>k_0$ so that also $A_k^{-1}\toto A^{-1}$ as $k\to\infty$.
At this point it is clear that $A^{-1}\in\BDO$ since
$A_k^{-1}\in\BDO$ by Proposition \ref{prop:BO} and since $\BDO$ is
closed.

For an explicit approximation of $A^{-1}$ by banded matrices, it
remains to look at banded approximations of $A_k^{-1}$ and to use
(\ref{eq:horiz}). Therefore, for every $k>k_0$, let $B_1^{(k)},
B_2^{(k)},...$ be the banded approximations of $A_k^{-1}$ from
Proposition \ref{prop:BO}, where $B_n^{(k)}\in\BO_{3nk}$ and
$\|A_k^{-1}-B_n^{(k)}\|$ is bounded by the right-hand side of
(\ref{eq:BO}) with $M,m$ and $\kappa$ replaced by $M_k:=\|A_k\|$,
$m_k:=1/\|A_k^{-1}\|$ and $\kappa_k:=M_k/m_k$, respectively.
Finally, put
\begin{equation} \label{eq:kappapm}
\kappa_k^\pm\ :=\ \frac{M\pm 2\delta_k}{m/\alpha_k}\ =\
\alpha_k\left(\kappa \pm \frac{2\delta_k}m\right)
\end{equation}
for every $k\in\N_0$. From $\|A-A_k\|\le 2\delta_k$ and
$\|A_k^{-1}\|\le \alpha_k \|A^{-1}\|$ we get $M-2\delta_k\le M_k\le
M+2\delta_k$ and $m_k\ge m/\alpha_k$. Hence
\begin{eqnarray}
\nonumber \|A_k^{-1}-B_n^{(k)}\|& \le&
\frac{M_k}{m_k^2}\left(\frac{M_k^2-m_k^2}{M_k^2+m_k^2}\right)^{n+1}
\ \le\
\frac{M+2\delta_k}{(m/\alpha_k)^2}\left(\frac{(M+2\delta_k)^2-(m/\alpha_k)^2}{(M-2\delta_k)^2+(m/\alpha_k)^2}\right)^{n+1}\\
\label{eq:vert}&=& \frac{(\kappa_k^+)^2}{M+2\delta_k}
\left(\frac{(\kappa_k^+)^2-1}{(\kappa_k^-)^2+1}\right)^{n+1}\ =\
\frac{\alpha_k^2\left(\kappa +
\frac{2\delta_k}m\right)^2}{M+2\delta_k}
\left(\frac{\alpha_k^2\left(\kappa +
\frac{2\delta_k}m\right)^2-1}{\alpha_k^2\left(\kappa -
\frac{2\delta_k}m\right)^2+1}\right)^{n+1}.
\end{eqnarray}
Bounding
$\|A^{-1}-B_n^{(k)}\|\le\|A^{-1}-A_k^{-1}\|+\|A_k^{-1}-B_n^{(k)}\|$
by (\ref{eq:horiz})+(\ref{eq:vert}) completes the proof.
\end{proof}

\section{The Wiener algebra $\W$ is inverse closed} \label{sec:Wiener}
The term 'Wiener algebra' is commonly used for the set $W(\T)$ of
all functions $f(t)=\sum_{n\in\Z}f_nt^n$ on the unit circle $\T$
whose sequence of Fourier coefficients $\hat f=(f_n)_{n\in\Z}$ is in
$\ell^1$ (two-sided infinite), equipped with pointwise addition and
multiplication and with the norm $\lW f\rW:=\sum |f_n|$. Norbert
Wiener's famous theorem says that if $f\in W(\T)$ is invertible as a
continuous function, i.e. $f$ vanishes nowhere on $\T$, then
$f^{-1}=1/f$ is in $W(\T)$ as well, showing that $W(\T)$ is inverse
closed.

In fact, Wiener's theorem is a special case of our Theorem
\ref{th:W} saying that 'our' Wiener algebra $\W$ is inverse closed!
It follows if we apply Theorem \ref{th:W} to a two-sided infinite
matrix with constant diagonals. To see this take $\I=\Z$ and
associate with every function $f\in W(\T)$ the so-called {\sl
Laurent matrix}
\begin{equation} \label{eq:Laurent}
L(f)\ =\ (f_{i-j})_{i,j\in\Z}\ =\ \left(\begin{array}{ccccc}
\ddots&\ddots&\ddots&\ddots&\ddots\\
\ddots&f_0&f_{-1}&f_{-2}&\ddots\\
\ddots&f_1&f_0&f_{-1}&\ddots\\
\ddots&f_2&f_1&f_0&\ddots\\
\ddots&\ddots&\ddots&\ddots&\ddots
\end{array}\right)
\end{equation}
which is sometimes also (wrongly) called a 'two-sided infinite
Toeplitz matrix'. Then, clearly,
\[
L(f)\in\W\qquad\textrm{and}\qquad \lW L(f)\rW\ =\ \sum_{n\in\Z}
|f_n|\ =\ \lW f\rW.
\]
Now for $g\in L^2(\T)$, let $\hat g=(g_n)_{n\in\Z}\in\ell^2$ denote
its sequence of Fourier coefficients and note that $L(f)\hat g=\hat
f\ast\hat g=\widehat{fg}$ acts as the operator of convolution by the
sequence $\hat f=(f_n)$. In other words, if $F:L^2(\T)\to\ell^2$ is
the Fourier transform $g\mapsto\hat g=(g_n)$ then
\[
L(f)\ =\ F\,M(f)\,F^{-1}
\]
with $M(f):L^2(\T)\to L^2(\T)$ denoting the multiplication operator
$g\mapsto fg$. The latter formula shows that $L(f)$ is invertible on
$\ell^2$ iff $M(f)$ is invertible on $L^2(\T)$ which is clearly the
case iff $f$ has no zeros on $\T$. In this case
$(L(f))^{-1}=FM(f^{-1})F^{-1}=L(f^{-1})$ holds, and Wiener's
statement, $f^{-1}\in W(\T)$, is equivalent to the inverse Laurent
matrix $(L(f))^{-1}=L(f^{-1})$ being in $\W$. Our Theorem \ref{th:W}
however says much more: For all matrices $A\in\W$, not just for
those with constant diagonals, the inverse $A^{-1}$, if it exists,
is in $\W$.

\begin{theorem} \label{th:W}
If $A\in\W$ is invertible then $A^{-1}\in\W$.
\end{theorem}
\begin{proof}
a) As in \S2, we start with the case $A\in\BO$, say $A\in\BO_k$ for
some $k\in\N_0$. Then, by (\ref{eq:BO}),
\[
\dist\big(\,A^{-1}\,,\,\BO_{3nk}\,\big)\ \le\ r_n\ :=\
\frac{\kappa^2}M\left(\frac{\kappa^2-1}{\kappa^2+1}\right)^{n+1}
\]
for every $n\in\N_0$, where $M=\|A\|$, $m=1/\|A^{-1}\|$ and
$\kappa=M/m=\|A\|\|A^{-1}\|$. Now, for every $j\in\Z$, let $d_j$
denote the supremum norm of the $j$-th diagonal of $A^{-1}$. By the
previous inequality and (\ref{eq:dist_d}) we get that
\[
\begin{array}{ccl}
d_j\le r_0&\textrm{for}& j\in\{\pm 1,...,\pm 3k\},\\
d_j\le r_1&\textrm{for}& j\in\{\pm (3k+1),...,\pm 6k\},\\
d_j\le r_2&\textrm{for}& j\in\{\pm (6k+1),...,\pm 9k\},\\
\vdots&&~\hspace{15mm}\vdots
\end{array}
\]
Summing up we have
\[
\lW A^{-1}\rW\ =\  \sum_{j\in\Z} d_j\ \le\
d_0+6kr_0+6kr_1+6kr_2+...\ =\ d_0+6k(r_0+r_1+r_2+...)\ <\ \infty
\]
since $r_n$ decays exponentially. So $A^{-1}\in\W$ if $A\in\BO$.

b) Now let $A\in\W$ be invertible and take $A_1, A_2, ...\in\BO$
such that $\lW A-A_i\rW \to 0$ as $i\to\infty$. Since
$(\W,\lW\cdot\rW)$ is a Banach algebra we know that for sufficiently
large $i$ also $A_i$ is invertible and $\lW A^{-1}-A_i^{-1}\rW\to 0$
as $i\to\infty$. Together with part a) and the closedness of
$(\W,\lW\cdot\rW)$ this proves the theorem.
\end{proof}
\bigskip

Note that one corollary of Theorem \ref{th:W} is that if, for some
fixed $p\in[1,\infty]$, an operator $A:\ell^p\to\ell^p$ with matrix
representation in $\W$ is invertible then its inverse is again given
by a matrix in $\W$ and $A^{-1}$ therefore (see (\ref{eq:WinBDO}))
acts boundedly on all spaces $\ell^p$ with $p\in[1,\infty]$. So
invertibility and spectrum of such operators $A$ do not depend on
the particular choice of $p$. In \cite{Li:Wiener} it is shown how
this result can be used to prove that even the property of being a
Fredholm operator (including the value of the Fredholm index) and
hence the essential spectrum of $A$ does not depend on
$p\in[1,\infty]$ if $A\in\W$.

\section{Some generalizations} \label{sec:gen}
\subsection{Generalized banding}
It is easy to see that our results generalize well beyond $\BDO$ and $\W$.
To see what we mean, let $\rh$ be a metric on $\I$ and define the set
$\BO^\rh$ of $\rh$-generalized banded operators as the set of all
$A=(a_{ij})_{i,j\in\I}\in \BL$ with $a_{ij}=0$ for all $i,j\in\I$ with
$\rh(i,j) > k$, for some fixed $k$. It is easy to see that $\BO^\rh$ is
also closed under addition and multiplication and taking adjoints.
Examples of interesting metrics $\rh$ other than $\rh(i,j)=|i-j|$ are
obtained by taking a sequence $(\bx_i)_{i\in\I}$ of pairwise different
elements from another metric space $(X,d)$ (e.g. $X=\R^n$ with the
Euclidean norm) and putting
\[
\rh(i,j)\ :=\ d(\bx_i,\bx_j),\qquad i,j\in\I.
\]
These generalizations are interesting in statistics as we shall see
in Section \ref{sec:stat}. If we let $\BDO^\rh$ be the closure of
$\BO^\rh$ in $\BL(\ell^2)$, Lemma \ref {lem:BOpos}, Proposition \ref{prop:BO}
and Theorem \ref{th:BDO} go over verbatim to $\BO^\rh$ and $\BDO^\rh$.

If we now define the generalized Wiener norm by
$$\lW A \rW_\rh \ := \ \sum_{k \in \Z} d_k^{(\rh)},$$
where $d_k^{(\rh)} = \sup \{ |a_{ij} |: i,j \in \I, \ \rh(i,j)=k \}$,
we can obtain an exact generalization of Theorem \ref{th:W}. An obvious
application of this result is the generalization of Wiener's theorem to
analytic functions
$$ f:\ \bt\quad\mapsto\quad \sum_{\bn \in \Z^p} f_{\bn}\bt^{\bn}$$
of several variables, where $\bt=(t_1,\ldots,t_p)\in\T^p$,
$\bn=(n_1,\ldots,n_p)\in\Z^p$ and $\bt^{\bn} = \prod^p_{j=1} t_j^{n_j}$.

\subsection{Banding up to a permutation}
Let $\pi$ be a permutation, that is, a 1-1 mapping of $\I$ onto itself
which is clearly representable as an operator
$(x_i)_{i\in\I}\mapsto (x_{\pi(i)})_{i\in\I}$ of norm 1 from $\ell^2$ to
$\ell^2$, which we also denote by $\pi$. Define
$$\BO^{(\pi)} = \{ A : \pi A \pi^* \in BO \},$$
where $\pi^*$ is the adjoint/inverse of $\pi$.
$\BO^{(\pi)}$ is closed under addition and multiplication but {\em not}
under taking adjoints since $A \in \BO^{(\pi)} \Longrightarrow A^* \in
\BO^{(\pi^*)}$. So we can't expect a generalization of Theorem \ref{th:BDO}.
We can, however, generalize Lemma \ref{lem:BOpos} and obtain a special
case of Theorem \ref{th:BDO}. Specifically, let
$$
\HBO^{(\pi)} \ :=\ \{ A  \in \BO^{(\pi)}: A \mbox{ is hermitian positive definite} \}
$$
and $\HBDO^{(\pi)}$ be its closure. Then Theorem \ref{th:BDO} generalizes to
$\HBDO^{(\pi)}$ if we replace $B^{(\pi)}_{3nk} $ by $B^{(\pi)}_{nk}$ (with
obvious notation changes).

We can also in an obvious way obtain the same conclusions for generalized
banding. More important from a statistical point of view is the following
generalization. Let
$$
\HBO^{\mbox{\tiny \sf PERM}} \ := \ \bigcup_{\pi} \HBO^{(\pi)}.
$$
Then $\HBO^{\mbox{\tiny \sf PERM}}$ is only closed for addition, scalar
multiplication, and taking powers. However, an examination of the proof
of Theorem \ref{th:BDO} shows that these properties are sufficient to
arrive at the same generalization for $\HBO^{\mbox{\tiny \sf PERM}}$
and its closure $\HBDO^{\mbox{\tiny \sf PERM}}$ as we did for
$\HBO^{(\pi)}$ and its closure. Again, everything carries over
verbatim to generalized banding. These results, particularly the last,
are of interest in statistics since, for reasons to become apparent,
it is desirable to define classes of covariance matrices such that
matrices and their inverses necessarily obey the same definition of sparseness.

\section{Applications to probability theory} \label{sec:Gauss}

\subsection{The closures of banded self-adjoint positive definite operators
and Gaussian processes}
A Gaussian process is a sequence of random variables,
$\{X_j:j\in\Z\}$, on a probability space whose finite
dimensional joint distributions are Gaussian. Without loss of
generality, we take $EX_j= 0$ for all $j$, so that the joint
distributions are determined by the matrices
$$
\Sigma_{m,n}\ :=\ E \X_m^n [ \X_m^n]^T\ =\ [ EX_iX_j ]_{i,j=m}^n, \qquad m,n \in\Z,
$$
where we put $\X_m^n := (X_m,\ldots,X_n)^T$. We extend these notations
to $m=-\infty$ and $n=\infty$ by introducing the two-sided infinite vector
$\X_{-\infty}^\infty := (\dots,X_{-1},X_0, X_{1},\dots)^T$ and matrix
$$
\Sigma\ :=\ E \X_{-\infty}^\infty [ \X_{-\infty}^\infty]^T\ =\ [ EX_iX_j ]_{i,j=-\infty}^\infty\ .
$$
For the rest of our discussion we assume that the infinite matrix
$\Sigma$ acts as a bounded operator from $\ell^2$ to $\ell^2$. A {\sl regular
process} is one for which
\[
E \in \mc B_{-\infty}\ :=\ \bigcap^{\infty}_{t=-\infty} \mc
B^{t}_{-\infty} \qquad \Longrightarrow\qquad P(E)\in\{0,1\},
\]
where $\mc B^n_m$ is the $\sigma$ field generated by $\X_m^n$ or
equivalently
$$\lim_{p \to \infty} \sup \big\{ | P(AB) -P(A)P(B)|: \ B \in
\mc B^{p}_{-\infty} \big\} = 0
$$
for all $A \in \mc B_{-\infty}$ and all $p$. A {\sl strongly mixing}
process is one for which
$$
\lim_{p \to \infty} \beta(p)=0,\quad \mbox{ where }\quad \beta(p)= \sup
\big\{ | P(AB) -P(A)P(B)|: \ A \in \mc B^m_{-\infty}, \ B \in\mc
B^{\infty}_{m+p} , \ m\in\Z \big\} \ .
$$
Let $P_{m,n}$ be the joint distribution of $\X_m^n$, the probability
measure induced on $\mc B_m^n$ by $P$. Let $P_m^{m+n}$ be the
regular conditional probability measure on $\mc B_{m+p}^\infty$
given $\mc B_{-\infty}^m$. If $\|\mu\|_{TV}$ denotes the total
variation of a signed measure $\mu$, let
\[
\beta(m,p)\ :=\ E \|P_m^{m+p}-P_{m+p,\infty}\|_{TV}.
\]
A {\sl beta mixing} (or {\sl absolutely regular}) process is one such that
\[
\lim_{p\to\infty} \sup_m \beta(m,p)\ =\ 0.
\]
As noted in \cite{IbragRoz}, beta mixing implies strong mixing. The
converse is not true as Example \ref{ex:strong_not_beta} below
shows. On the other hand, let $Q_{m,n}$ be the probability
distribution induced by $P$ on the $\sigma$ field generated by $\mc
B_{-\infty}^m$ and $\mc B_n^\infty$. Let $\overline Q_{m,n}$ be the
product probability on the same $\sigma$ field with marginals
$P_{-\infty,m}$, $P_{n,\infty}$. Then Lemma 2 of \cite[p.
118]{IbragRoz} states
\begin{equation} \label{eq:T1}
\beta(m,p)\ =\ \frac 12\ \|Q_{m,m+p}-\overline Q_{m,m+p}\|_{TV}.
\end{equation}


A mean 0 process is {\sl linearly regular} if
$$E(X_{p+m+1}| X_{-\infty}^m)\ \to\ 0 $$
as $p \to \infty$, uniformly in $m$. A detailed discussion of these
concepts is in Ibragimov and Rozanov \cite{IbragRoz} primarily in
the context of stationary processes.

For Gaussian processes, linear regularity and regularity are
equivalent, see \cite[p.112]{IbragRoz}. We connect with our previous
results via

\begin{theorem} \label{th:51}
If $\Sigma$ has a bounded inverse and belongs to $\BDO$ then
$\X_{-\infty}^{\infty}$ is regular, and so is the process
corresponding to $\Sigma^{-1}$.
\end{theorem}

\begin{proof}
Recall $\X^b_a = (X_a,X_{a+1}, \ldots, X_b)^T$ and let
$$
\Sigma (a,b) \ = \ E\X^b_a [ \X^b_a]^T\ =\ [ EX_iX_j ]_{i,j=a}^b
$$
denote the $(b-a+1)\times (b-a+1)$ diagonal block of $\Sigma$.
Moreover let
$$\bsig (a,b,c) \ =\ E \X^b_a  X_c \ .$$
Then
\begin{eqnarray}
\nonumber E(X_{m+p+1}| \X^n_{-m}) &=& [\X^n_{-m}]^T \Sigma^{-1}(-m,n) \bsig(-m,n,n+p+1)\\
\label{eq6.1} E\big[ E(X_{n+p+1}| \X^n_{-m}) \big]^2 &=&
\bsig^T(-m,n,n+p+1)\Sigma^{-1}(-m,n) \bsig(-m,n,n+p+1) \ .
\end{eqnarray}
Since $\Sigma$ is a member of $\BDO$ we can find $B_{\ve}$ banded of
width $k(\ve)$ such that $\| B_{\ve} - \Sigma \| \leq \ve $. Let
$\big(X_{-m}(\ve),\ldots, X_{n+p+1}(\ve), \ldots \big)$ be a
Gaussian process with covariance operator $B_{\ve}$, (a moving
average process). Then,
\begin{eqnarray}
\nonumber \| B_{\ve} (-m,n) - \Sigma (-m,n) \| &\leq& \ve \\
\label{eq6.2} | \bsig_{\ve} (-m,n,n+p+1) - \bsig(-m,n,n+p+1) |
&\leq& \ve,
\end{eqnarray}
where $\bsig_{\ve}(a,b,c)= E\X^b_a  X_c (\ve)$ and $|\cdot|$ is the
Euclidean norm. By construction, \be \label{eq6.3}
\bsig_{\ve}(-m,n,n+p+1)  = \0 \quad \mbox{ if } \quad p \geq k(\ve)
\ . \ee Then, from \eqref{eq6.2} and \eqref{eq6.3},
\begin{equation}
\label{eq6.4} |\bsig^T(-m,n,n+p+1)\Sigma^{-1}(-m,n)\bsig(-m,n,n+p+1)
|\ \le\ \ve \| \Sigma^{-1} \| \ve
\end{equation}
for all $p \geq k(\ve)$.
The main result follows from \eqref{eq6.1} and \eqref{eq6.4}. The
corresponding statement for $\Sigma^{-1}$ is a consequence of Theorem \ref{th:BDO}.
\end{proof}

\subsection{Beta mixing and Frobenius closure of banded operators}
In general, $\Sigma \in \BDO$ having a bounded inverse does not
imply beta mixing (see Example \ref{ex:strong_not_beta} below).
But below we prove that beta mixing is equivalent to a condition
on the off-diagonal decay of $\Sigma$ which can be related to the
closure of $\BO$ in a type of Frobenius norm.

In the negative, there are results of Kolmogorov and Rozanov \cite{KolmRoz} for
symmetric positive definite Toeplitz matrices $\Sigma$, i.e. one-sided infinite
versions of \eqref{eq:Laurent} with $f_n=f_{-n}$, showing that strong (and hence beta)
mixing does not hold if the associated symbol $f: t\in\T \mapsto \sum_{n\in\Z}f_nt^n\in\R$
has discontinuities of the first kind. (Recall from Section \ref{sec:Wiener}
that $\Sigma$ is: bounded iff $f\in L^\infty(\T)$, invertible iff $1/f\in L^\infty(\T)$,
and it is in $\BDO$ iff $f$ is continuous.)


On the positive side, Ibragimov and Rozanov \cite[p.129]{IbragRoz}
establish

\begin{theorem} \label{th:52}
If $f$ is the symbol of a (one- or two-sided infinite) Toeplitz matrix
$\Sigma$ corresponding to a stationary Gaussian process $X$ and
\[
f(e^{ix})\ =\ |P(e^{ix})|^2\,a(x),\qquad x\in (-\pi,\pi),
\]
where $P$ is a polynomial with zeros, if any, only on the unit circle
and $\log a(\cdot)$ belongs to the Sobolev space $W^{\frac 12,2}$
then $X$ is beta mixing and conversely.
\end{theorem}

Recall that
\[
W^{s,2}\ =\ \left\{ b(x) = \sum_{k\in\Z} a_k e^{ikx}\ :
\ \sum_{k\in\Z} |k|^{2s}\,|a_k|^2 < \infty\right\}.
\]
Note that if $f$ is bounded above and away from zero then one can take
$P\equiv 1$, and the condition $\log a(\cdot)\in W^{\frac 12,2}$
is equivalent to $a\in W^{\frac 12,2}$ and $a$ bounded away from zero.
To see the latter, note that $W^{\frac 12,2}$ can be equivalently
characterized by the Sobolev-Slobodeckij norm (e.g. \cite{RunstSickel}),
in which it becomes clear that with $f$ also powers of $f$ and hence,
by closedness, also $\log f$ (if $f$ is bounded above and away from zero)
and $\exp f$ are in $W^{\frac 12,2}$.

Here is an example of a strong but not beta mixing stationary process.
\begin{example} \label{ex:strong_not_beta}
For $k\in\Z$, let $a_k=1/\sqrt{|k|}$ if $|k|\in \{1^4, 2^4, 3^4, ... \}$,
and $a_k = 0$ otherwise. Now look at the symbol function
$f(t) = \sum_{k\in\Z} a_k t^k$ defined on the unit circle $\T$.
Because of $a_k = a_{-k}$, the symbol $f$ is real-valued. One moreover has
\[
\sum_{k\in\Z} |a_k|\ =\ 2 \sum_{m\in\N} \frac 1{\sqrt{m^4}}
\ =\ 2 \sum_{m\in\N} \frac 1{m^2}\ =\ \frac{\pi^2}3\ <\ \infty,
\]
so that $f$ is in the Wiener class $W(\T)$ that we discuss at the
beginning of Section \ref{sec:Wiener}, and hence $f$ is continuous.
In particular, $f$ is bounded with
\[
|f(t)|\ =\ \left| \sum_{k\in\Z} a_k t^k \right|\ \le\ \sum_{k\in\Z} |a_k|
\ =\ \frac{\pi^2}3 \ <\  4,\qquad t \in \T,
\]
so that $f(t)\in(-\frac{\pi^2}3,\frac{\pi^2}3)\subset (-4,4)$ for all $t\in\T$.
However, $f$ is not in the Sobolev space $W^{\frac 12 , 2}$ since
\[
\sum_{k\in\Z} |k| |a_k|^2 \ =\ 2 \sum_{m\in\N} m^4 \frac 1{m^4}\ =\ \infty.
\]
Putting $g(t):=f(t)+4$, we get that the Fourier coefficients $b_k$ of
$g$ coincide with $a_k$, except $b_0$ which is $4$. Now the range
of $g$ is in $(4-\frac{\pi^2}3,4+\frac{\pi^2}3)\subset (0,8)$, so that
$g$ is positive, bounded away from zero, and continuous, whence the
associated process with covariance matrix $\Sigma = (b_{i-j})_{i,j}$
is strong mixing. But the process is not beta mixing since
$g\not\in W^{\frac 12 , 2}$.
\end{example}

We now give a generalization of Theorem \ref{th:52} to arbitrary
bounded covariance matrices $\Sigma$. Implicitly, the result is
essentially in Lemmas 2-5 of \cite[\S IV.4]{IbragRoz} but we give a
full statement and proof here for completeness. We denote the
entries of our infinite covariance matrix $\Sigma$ by
$\sigma_{ij}=EX_iX_j$ for $i,j\in\Z$.

\begin{theorem} \label{thm:betamixing}
Suppose $\Sigma\in\BL$ is invertible and $E\X_{-\infty}^{\infty} =
\0$. Then $\X^{\infty}_{-\infty}$ Gaussian is beta mixing iff
\begin{equation} \label{neweq4.6}
\sup_{n\in\Z}\ \sum^n_{i=-\infty} \sum^{\infty}_{j=p+n}
\sigma_{ij}^2\ \to\ 0  \quad \textrm{as} \quad p\to\infty.
\end{equation}
\end{theorem}

A simpler sufficient condition is given by
\begin{equation} \label{eq6.7}
\gamma(p)<\infty\quad \textrm{for some} \quad p\ge
1,\qquad\textrm{where}\qquad \gamma(p)\ :=\
\sum^{\infty}_{i=-\infty} \sum^{\infty}_{j=p+i} \sigma_{ij}^2,
\end{equation}
since \eqref{eq6.7} $\Rightarrow$ \eqref{neweq4.6}.

\medskip
\begin{proof}
Note that with $\Sigma$ also $\Sigma^{-1}$ is in $\BL$ and
put $M:=\max(\|\Sigma\|,\|\Sigma^{-1}\|)<\infty$.
By \eqref{eq:T1}, to prove the theorem we need only show that
\eqref{neweq4.6} is equivalent to
\begin{equation} \label{eq:T2}
\sup_m \|Q_{m,m+p}-\overline Q_{m,m+p}\|_{TV}\ \to\ 0\qquad\textrm{as}\qquad p\to\infty.
\end{equation}
Moreover, it is easy to see that \eqref{eq:T2} is equivalent to
\begin{equation} \label{eq:T3}
\sup_{m,n,k}\|Q_{m,n,p,k}-\overline Q_{m,n,p,k}\|_{TV}\ \to\
0\qquad\textrm{as}\qquad p\to\infty,
\end{equation}
where $Q_{m,n,p,k}$ and $\overline Q_{m,n,p,k}$ are the
distributions $Q_{m,m+p}$ and $\overline Q_{m,m+p}$ restricted to
$\mc B_m^n$, $\mc B_{n+p}^{n+p+k}$.

We consider $\X^{(1)}:=(X_m,\ldots,X_n)^T$, $\X^{(2)} :=
(X_{n+p+1},\ldots,X_{n+p+k})^T$. Let $f$ denote the joint density of
$(\X^{(1)}, \X^{(2)})$ so that $f$ corresponds to
$$
S = \left(
\begin{array}{cc}
\Sigma_{11} & \Sigma_{12} \\
\Sigma_{21} & \Sigma_{22}
\end{array}
\right),
$$
the covariance matrix of $(\X^{(1)}, \X^{(2)})$ blocked out. Thus,
$$
\X^{(1)} \sim N_{n-m+1}(0,\Sigma_{11})\qquad \textrm{and} \qquad
\X_2 \sim N_k(0,\Sigma_{22}).
$$
Let $g$ correspond to
$$S_0  := \left(
\begin{array}{cc}
\Sigma_{11} & 0 \\
0  & \Sigma_{22}
\end{array}
\right), $$
the covariance matrix of $(\X^{(1)}, \X^{(2)})$ if $\X^{(1)}$ and $\X^{(2)}$ are independent.
Let $P_f, P_g$ be the probability distributions of $\X^{(1)}$ and $\X^{(2)}$ and
$$\| P_f - P_g \|_{TV} \ :=\ \frac{1}{2} \int |f -g|$$
be the variational norm. We suppress the dependence of $f,g$ on $m,p,k$ in what follows.
Let
$$H^2(P_f,P_g)\ :=\ 2\left( 1 - \int [fg]^{\frac{1}{2}} \right) \ =:\
2\big(1 - A(f,g) \big)  $$
be the (squared) Hellinger metric. It is well known that
\[  
\frac 12 H^2(P_f,P_g)\ \le\ \|P_f-P_g\|_{TV}\ \le\ \sqrt 2 H^2(P_f,P_g).
\]  
Thus we can replace $\|Q_{m,n,p,k}-\overline Q_{m,n,p,k}\|_{TV}$ by
$H^2(Q_{m,n,p,k},\overline Q_{m,n,p,k})$ in \eqref{eq:T3}.

Our argument will consist of bounding $\frac 12 H^2 = 1-A$ above and
below by functions $\underline a_{m,n,p,k}$ and $\overline
a_{m,n,p,k}$ of $\{\sigma_{ij}^2 : m\le i\le n, n+p+1\le j\le
n+p+k\}$ (see \eqref{eq:squeeze} below) such that
\[
\sup_{m,n,k} \underline {\overline a}_{m,n,p,k}\ \to\
0\qquad\textrm{as}\qquad p\to\infty
\]
iff \eqref{neweq4.6} holds. To do this, we have to compute $A(f,g)$.



Note that $\| \cdot \|_{TV}$ and $H^2(\cdot,\cdot)$ are invariant
under regular linear transformations $\X^{(1)} \to T_1\X^{(1)}$ and
$\X^{(2)} \to T_2\X^{(2)}$. For the choice of these matrices $T_1$
and $T_2$, suppose, by the spectral theorem, that $\Sigma_{jj} =
Q_j^T \Lambda_j Q_j$, $j=1,2$, where $Q_j$ are orthogonal and
$\Lambda_j$ are diagonal, and put
\[
T_1 \ =\  \Lambda_1^{-\frac{1}{2}} Q_1 \qquad\textrm{and}\qquad
T_2 \ =\ \Lambda_2^{-\frac{1}{2}} Q_2.
\]
Replacing $(\X^{(1)}, \X^{(2)})$ by $(T_1\X^{(1)}, T_2\X^{(2)})$ in
the above, we get that $\Sigma_{11} = I_1$ and $ \Sigma_{22} = I_2$
are the $(n-m+1)\times (n-m+1)$ and $k \times k$ identity,
respectively. We shall establish the theorem in this case and then
derive the general case.

If $\Sigma_{11} = I_1$ and $\Sigma_{22} = I_2$ then, in
corresponding block notation,
$$S^{-1} = \left(
\begin{array}{cc}
\Sigma^{11} & \Sigma^{12} \\
\Sigma^{21} & \Sigma^{22}
\end{array}   \right)
$$
with
\begin{eqnarray}
\nonumber \Sigma^{11} &=& (I_1 - \Sigma_{12} \Sigma_{21} )^{-1}\\
\nonumber \Sigma^{22} &=& (I_2 - \Sigma_{21} \Sigma_{12} )^{-1}\\
\label{eqIII} \Sigma^{12} &=& - (I_1 - \Sigma_{12} \Sigma_{21} )^{-1} \Sigma_{12}\\
\nonumber \Sigma^{21} &=& - (I_2 - \Sigma_{21} \Sigma_{12} )^{-1} \Sigma_{21}
\end{eqnarray}
and the determinant of $S$ is equal to
\begin{equation} \label{eqV}
|S|\ =\ | I_1 -\Sigma_{12} \Sigma_ {21} |\ =\ | I_2 -\Sigma_{21}
\Sigma_ {12} |
\end{equation}
since $\Sigma_{12}\Sigma_{21}$ and $\Sigma_{21} \Sigma_{12}$ have
the same nonzero eigenvalues.
It holds that
\begin{equation}
\label{neweq4.78} A(f,g) \ =\ \int (fg)^{\frac{1}{2}} \ =\ |S|^{-
\frac{1}{4}}\ (2\pi)^{-\frac{n-m+1+k}2} \int e^{-\frac 12\bx^T\frac
12(S^{-1}+S_0^{-1})\bx} d\bx \ =\ \ |S|^{-\frac{1}{4}}\ \left|\frac
12(S^{-1}+S_0^{-1})\right|^{-\frac{1}{2}}
\end{equation}
with $\bx={\bx_1\choose \bx_2}\in\R^m\times\R^k$ and
\begin{eqnarray}
\nonumber S^{-1}+S_0^{-1} &=&  \left(\begin{array}{cc}
\Sigma^{11} + I_1 & \Sigma^{12}\\
\Sigma^{21} & \Sigma^{22} + I_2
\end{array}\right)\\
\label{eq:M1122} & =& \left(\begin{array}{cc}
(I_1 - \Sigma_{12} \Sigma_{21} )^{-1} + I_1 & - (I_1 - \Sigma_{12} \Sigma_{21} )^{-1} \Sigma_{12}\\
- (I_2 - \Sigma_{21} \Sigma_{12} )^{-1} \Sigma_{21} & (I_2 -
\Sigma_{21} \Sigma_{12} )^{-1} + I_2
\end{array}\right)\ =:\
\left(\begin{array}{cc}
M_{11} & M_{12}\\
M_{21} & M_{22}
\end{array}\right).
\end{eqnarray}
To compute \eqref{neweq4.78}, recall \eqref{eqV}, \eqref{eq:M1122}
and the standard formula for block determinants
\begin{equation} \label{eq:blockdet}
|S^{-1}+S_0^{-1}| = |M_{11}|\ |M_{22}-M_{21}M_{11}^{-1}M_{12}|.
\end{equation}
W.l.o.g suppose $n-m+1\le k$. Now let $\lambda_m,\ldots,\lambda_n$
be the eigenvalues of $\Sigma_{12} \Sigma_{21}$. Then $\Sigma_{21}
\Sigma_{12}$ has the same $n-m+1$ eigenvalues and the rest are
zeros. Let $x$ be an eigenvector of $\Sigma_{21} \Sigma_{12}$
corresponding to $\lambda$. Then $\Sigma_{12} x$ is an eigenvector
of $\Sigma_{12} \Sigma_{21}$ corresponding to the same $\lambda$.
Consequently, for such an $x$, we have
\begin{eqnarray*}
M_{22}x &=& ((I_2 - \Sigma_{21} \Sigma_{12} )^{-1} + I_2)x \ =\
((1-\lambda)^{-1}+1)x\ =\ (1-\lambda)^{-1}(2-\lambda)x
\end{eqnarray*}
and
\begin{eqnarray*}
M_{21}M_{11}^{-1}M_{12}x &=& (I_2 - \Sigma_{21} \Sigma_{12} )^{-1}
\Sigma_{21}((I_1 - \Sigma_{12} \Sigma_{21} )^{-1} + I_1)^{-1}(I_1 -
\Sigma_{12} \Sigma_{21} )^{-1} \Sigma_{12}x\\
&=& (I_2 - \Sigma_{21} \Sigma_{12} )^{-1} \Sigma_{21}((1 -
\lambda)^{-1} + 1)^{-1}(1 -
\lambda )^{-1} \Sigma_{12}x\\
&=& (1 - \lambda )^{-1} ((1 - \lambda)^{-1} + 1)^{-1}(1 -
\lambda )^{-1} \lambda x\\
&=& (1 - \lambda )^{-1} (2 - \lambda )^{-1} \lambda x.
\end{eqnarray*}
Taking this together with \eqref{eq:M1122} and \eqref{eq:blockdet},
we get
\begin{eqnarray}
\nonumber \left|\frac 12(S^{-1}+S_0^{-1})\right| &=& \prod_{j=m}^n
\frac 12\left(\frac 1{1-\lambda_j}+1\right)\ \prod_{j=1}^k \frac
12\left(\frac{2-\lambda_j}{1-\lambda_j}-\frac{\lambda_j}{(1-\lambda_j)(2-\lambda_j)}\right)\\
\label{eq:detM} &=& \prod_{j=m}^n \frac 14\ \frac
{2-\lambda_j}{1-\lambda_j}\
\frac{(2-\lambda_j)^2-\lambda_j}{(1-\lambda_j)(2-\lambda_j)} \ =\
\prod_{j=m}^n \frac {(4-\lambda_j)(1-\lambda_j)}{4(1-\lambda_j)^2}\
=\ \prod_{j=m}^n \frac {1-\frac{\lambda_j}4}{1-\lambda_j}\ ,
\end{eqnarray}
so that, by \eqref{eqV}, \eqref{neweq4.78} and \eqref{eq:detM},
\begin{equation} \label{eq:Afg_final}
\nonumber A(f,g) \ =\ |S|^{-\frac{1}{4}}\ \left|\frac
12(S^{-1}+S_0^{-1})\right|^{-\frac{1}{2}}\ =\ \left(\prod_{j=m}^n
(1-\lambda_j)\right)^{-\frac 14}\left(\prod_{j=m}^n \frac
{1-\frac{\lambda_j}4}{1-\lambda_j}\right)^{-\frac 12} \ =\
\prod_{j=m}^n
\frac{(1-\lambda_j)^{\frac14}}{(1-\frac{\lambda_j}4)^{\frac12}} \ ,
\end{equation}
where we recall that $0\le \lambda_j < 1$ for all $j$ since
$\Sigma_{12} \Sigma_{21}$ is positive semi-definite and $(\Sigma^{11})^{-1} = I_1 -
\Sigma_{12} \Sigma_{21}$ is positive definite.

Now we put
$$t\ :=\ {\rm Trace}(\Sigma_{12}\Sigma_{21})\ =\ \sum_{j=m}^n \lambda_{j}$$
and note that our condition \eqref{neweq4.6} is equivalent to
\[
t\ =\ {\rm Trace}(\Sigma_{12}\Sigma_{21}) \ =\ \sum^{n}_{i=m} \left[
\sum^{n+p+k}_{j=n+p+1} \sigma_{ij}\, \sigma_{jl} \right]_{ i=l} =
\sum^{n}_{i=m} \sum^{n+p+k}_{j=n+p+1} \sigma_{ij}^2 \ \to\ 0
\quad\mbox{ as }\quad p \to \infty
\]
uniformly in $m,n,k$. The inequalities
\[ 
\prod_{j=m}^n (1-\lambda_j)\ \ge\ 1\,-\,\sum_{j=m}^n
\lambda_j,\qquad 0\le \lambda_m,...,\lambda_n\le 1
\] 
(as can be seen by induction over the number of terms) and
\[ 
\frac{1-\lambda}{(1-\frac{\lambda}4)^2}\ \le\ e^{-\frac 12 \lambda},\qquad \lambda\ge 0
\] 
(which is easily checked using basic calculus), together with \eqref{eq:Afg_final}, yield
\begin{equation} \label{eq:squeeze}
(1-t)^{\frac 14}\ \le\ \prod_{j=m}^n (1-\lambda_j)^{\frac14}\ \le\
A(f,g)\ \le\ e^{-\frac 18 t}\ \le 1 .
\end{equation}
From \eqref{eq:squeeze} we get that $t\to 0$ implies $A(f,g)\to 1$. Conversely, by the right half of \eqref{eq:squeeze}, $A(f,g)\to 1$ implies $t\to 0$.

Thus the result is proved if $\Sigma_{11}$ and $\Sigma_{22}$ are the identity.

\medskip
\noindent {\bf General case. }
By the spectral theorem we noted we can find $Q_1$ and $Q_2$ orthogonal
such that
\[
\Sigma_{11} \ =\ Q_1^T\Lambda_1 Q_1\qquad\textrm{and}\qquad
\Sigma_{22} \ =\ Q_2^T\Lambda_2 Q_2,
\]
where $\Lambda_{1}$ and $\Lambda_{2}$ are diagonal. The transformation
of $\R^{m+k}$ by
$\left( \begin{array}{cc} Q_1 & 0 \\ 0 & Q_2 \end{array} \right)$
doesn't change Hellinger or variational distances and sends
$\Sigma_{ij} \mapsto \tilde{\Sigma}_{ij}$, where
\[
\tilde\Sigma_{11}=\Lambda_1,\quad
\tilde\Sigma_{22}=\Lambda_2,\quad
\tilde\Sigma_{12}= Q_1 \Sigma_{12} Q_2^T, \quad
\tilde\Sigma_{21}= Q_2 \Sigma_{21} Q_1^T \quad
\textrm{and hence}\quad
\tilde{\Sigma}_{12}\tilde{\Sigma}_{21}
= Q_1 \Sigma_{12} \Sigma_{21} Q_1^T.
\]
But then
\be
\label{neweq4.16}
\mbox{Trace}(\Sigma_{12} \Sigma_{21}) \ = \ \mbox{Trace}(\tilde{\Sigma}_{12}
\tilde{\Sigma}_{21} )
\ee
since eigenvalues are unchanged.
Sending
\[
\X^{(1)} \ \mapsto\ \Lambda_1^{-\frac{1}{2}} Q_1\X^{(1)}
\qquad\textrm{and}\qquad
\X^{(2)} \ \mapsto\ \Lambda_2^{-\frac{1}{2}} Q_2\X^{(2)}
\]
again doesn't change Hellinger and TV distances. If the resulting
covariance matrices are $\Sigma^*_{ij}$ then
\[
\Sigma^*_{11}\ =\ \Sigma^*_{22} \ =\ J
\qquad\textrm{and}\qquad
\Sigma^*_{12} \Sigma^*_{21} \ =\ \Lambda_1^{-\frac{1}{2}} \tilde{\Sigma}_{12}
\Lambda_2^{-1} \tilde{\Sigma}_{21} \Lambda_1^{-\frac{1}{2}}.
\]
Then
\begin{equation}
\label{neweq4.17} \mbox{Trace}(\Sigma^*_{12} \Sigma^*_{21})\ =\
\sum_{i=m}^n\sum_{j=n+p+1}^{n+p+k} \frac{\tilde{\sigma}_{ij}^2}{
\lambda^{(i)}_1 \lambda^{(i)}_2}
\end{equation}
where $\tilde{\Sigma}_{12} = ( \tilde{\sigma}_{ij} )$ and
$\lambda^{(i)}_1$, $\lambda^{(i)}_2$ are the diagonal elements of
$\Lambda_1,\Lambda_2$. But since
\[
M^{-1}\ \leq\ \|\Sigma^{-1}\|^{-1}\ \leq\ \lambda^{(i)}_1,
\lambda^{(i)}_2 \ \leq\ \|\Sigma\|\ \le \ M,\qquad i=m,...,n
\]
it follows that
$$M^{-2}  \ \mbox{Trace}(\Sigma_{12} \Sigma_{21}) \ \leq\ \mbox{Trace}
(\Sigma^*_{12} \Sigma^*_{21}) \ \leq \ M^2 \ \mbox{Trace}( \Sigma_{12}
\Sigma_{21}) $$
by \eqref{neweq4.16} and \eqref{neweq4.17}. 
Now the theorem follows.
\end{proof}

\begin{remark}
If we suppose $\Sigma$ to be a Laurent matrix, i.e.
$\Sigma=(a_{i-j})_{i,j\in\Z}$ with $a_k=a_{-k}$ for all $k$, then it
is easy to see that \eqref{neweq4.6} is equivalent to
\[ 
\sum^{\infty}_{i=1} \sum^{\infty}_{k=p+i} a_k^2\ \to\ 0,
\qquad\textrm{which evidently holds iff}\qquad
\sum^{\infty}_{i=1} \sum^{\infty}_{k=i+1} a_k^2\ <\ \infty,
\]
i.e. iff
\[ 
\label{eq:Sobolev}
\sum^{\infty}_{k=1} k\, a_k^2\ <\ \infty.
\] 
This is the $W^{\frac 12,2}$ condition from \cite{IbragRoz} (see
Theorem \ref{th:52}).
\end{remark}

We now prove a closure under inversion result similar to those of Sections
\ref{sec:BDO} and \ref{sec:Wiener}. Therefore, consider the cone of bounded
self-adjoint positive definite operators $\Sigma = (\sigma_{ij})$ with
bounded inverses and, for $m=0,1,2,...$, denote the subcone of all such
operators with
$$
\sum_{|i-j|\ge m} \sigma^2_{ij}\ <\ \infty
$$
by $F_m$. Define on $F_m$ an equivalence relation by
$$( \sigma_{ij} ) \equiv (\tau_{ij}) $$
iff $\sigma_{ij}\, \chi( |i-j| \geq m) = \tau_{ij}\, \chi( |i-j| \geq m)$,
where
\[
\chi(|i-j| \geq m)\ =\ \left\{
\begin{array}{cl} 1&\textrm{if } |i-j| \geq m,\\ 0&\textrm{otherwise.}
\end{array}\right.
\]
Let $\Sigma_0 := \big( \sigma_{ij}\ \chi( |i-j| \geq m) \big)$
correspond to such an equivalence class and define
$$
\| \Sigma_0 \|^2_{F_m} \ :=\ \sum_{|i-j| \geq m} \sigma^2_{ij}.
$$
The quotient cone, which we again denote by $F_m$, is closed under convex
combination and positive scaling but not under multiplication.
However it is closed under Schur multiplication
$$\Sigma_0 \ast \Ta_0 \ :=\ \big( \sigma_{ij}\, \tau_{ij}\, \chi( |i-j| \geq m)\big),$$
where $\Sigma_0, \Ta_0$ and $\Sigma_0 \ast \Ta_0$ are equivalence classes.
Clearly
\[ 
\| \Sigma_0 \ast \Ta_0 \|_{F_m} \ \leq\ \| \Sigma_0\|_{F_m} \| \Ta_0 \|_{F_m}.
\]  
Now note that $F_0 \subset F_1 \subset \cdots$ and put
$$F \ :=\ \bigcup_{m\ge 0} F_m,$$
so that a matrix $\Sigma=(\sigma_{ij})$ is in $F$ iff
$\|\Sigma\|_{F_m}<\infty$ for some $m\ge 0$. For
$\Sigma\in\cap_{m\ge 0} F_m = F_0$, one has that
\[
\sup_{m\ge 0} \|\Sigma\|_{F_m}^2\ =\ \|\Sigma\|_{F_0}^2\ =\
\sum_{i,j} \sigma_{ij}^2
\]
is the usual Frobenius norm.

We note the connection to $\BO$. Evidently, $\BO_m\subset F_{m+1}$
for all $m\ge 0$ and hence $\BO\subset F$. In fact, $F$ is the
closure of $\BO$ in the Frobenius metric $d_F(\Sigma,\Ta)\ :=\
\|\Sigma - \Ta\|_{F_0}$. To see the latter, note that $F$ is closed
under $d_F(\cdot,\cdot)$ and take $\Sigma=(\sigma_{ij})\in F$, so
that $\Sigma\in F_{m_0}$, i.e. $\sum_{|i-j|\ge m_0} \sigma^2_{ij} <
\infty$, for some $m_0\ge 0$. Then, putting $B_m:=\big(\sigma_{ij}\,
\chi( |i-j| < m)\big)\in\BO_m$ for $m=0,1,...$, it holds that
\[
d_F(\Sigma, B_m)^2\ =\ \|\Sigma-B_m\|_{F_0}^2\ =\
\|\Sigma-B_m\|_{F_m}^2\ =\ \sum_{|i-j|\ge m} \sigma_{ij}^2\ \to\
0\qquad\textrm{as}\qquad m\to\infty.
\]
We will see now that, besides $\BDO$ and $\W$ (see Theorem
\ref{th:BDO} and \ref{th:W}), which are the closures of $\BO$ under
$\|\cdot\|$ and $\lW\cdot\rW$, respectively, also the Frobenius
closure $F$ of $\BO$ is inverse closed.

\begin{theorem}
$F$ is closed under inversion.
\end{theorem}
As we have noted in Theorem \ref{thm:betamixing}, all operators in
$F$ are beta mixing so that we conclude that at least on $F$ beta
mixing is preserved under inversion.

\begin{proof}
Let $\Y_{-\infty}^\infty=(\ldots, Y_m,\ldots, Y_n,\ldots)^T$ have
covariance operator $\Sigma ^{-1}$ and, following our previous
practice from the proof of Theorem \ref{thm:betamixing}, consider
$\Y^{(1)} := (Y_m,\ldots,Y_n)^T$ and $\Y^{(2)} := (Y_{n+p+1},
\ldots,Y_{n+p+k} )^T$ and the corresponding covariance matrices,
$\Sigma^{11}$, $\Sigma^{12}$, $\Sigma^{21}$, $\Sigma^{22}$.
Following our previous argument we need only check that
$$ \mbox{Trace}(\Sigma^{12}\Sigma^{21}) \to 0$$
as $p \to \infty$ uniformly in $m,n,k$. As before we can reduce to
the case where $\Sigma^{11}$ and $\Sigma^{22}$ are the identity.
Hence, by formula \eqref{eqIII},
\[
\Sigma^{12} \ =\ -(I_1 - \Sigma_{12} \Sigma_{21} )^{-1}\Sigma_{12}
\qquad\textrm{and}\qquad
\Sigma^{21} \ =\ (\Sigma^{12})^T\ =\ -\Sigma_{21}(I_1 - \Sigma_{12} \Sigma_{21} )^{-1},
\]
so that
\be \label{eq6.20} \mbox{Trace}(\Sigma^{12}\Sigma^{21})\ =\
\mbox{Trace}(\Sigma^{21}\Sigma^{12})\ =\ \mbox{Trace} \big(
\Sigma_{21} (I_1 - \Sigma_{12} \Sigma_{21} )^{-2} \Sigma_{12} \big).
\ee Next, note that \be \label{eq6.22} \| \Sigma_{12} \Sigma_{21}
\|^2_F \ = \ \sum^n_{j=m} \lambda^2_j \ \leq\ \Big(\sum^n_{j=m}
\lambda_j \Big)^2 \ =\
\big(\,\mbox{Trace}(\Sigma_{12}\Sigma_{21})\,\big)^2, \ee where
$\lambda_m,\ldots,\lambda_n$ are the eigenvalues of
$\Sigma_{12}\Sigma_{21}$. But \be \label{eq6.23}
\mbox{Trace}(\Sigma_{12}\Sigma_{21}) \ = \ \sum^n_{j=m} \lambda_j \
=\ \sum^n_{i=m} \sum^{n+p+k}_{j=n+p+1} \sigma_{ij}^2 \ \to\ 0 \ee as
$p \to \infty$ since $\Sigma \in F_n$ for some $n$. Hence, for $p$
large enough, we can write \be \label{eq6.24} ( I_1 - \Sigma_{12}
\Sigma_{21} )^{-2} \ =\ \sum^{\infty}_{k=0} {-2 \choose k} (-
\Sigma_{12} \Sigma_{21} )^k \ =\ I_1 + 2\Sigma_{12} \Sigma_{21} + O
\big(\| \Sigma_{12} \Sigma_{21} \|^2_F \big) \ . \ee Then, by
\eqref{eq6.20} and \eqref{eq6.24},
$$
\mbox{Trace} (\Sigma^{12}\Sigma^{21}) \ =\ \mbox{Trace} (\Sigma_{12} \Sigma_{21} ) + 2\,\mbox{Trace} (\Sigma_{12}
\Sigma_{21} )^2 + O\big( \|\Sigma_{12} \Sigma_{21} \|^3_F \big).
$$
In view of \eqref{eq6.22} and \eqref{eq6.23}, the proof is complete.
\end{proof}

\section{Applications to Statistics} \label{sec:stat} This paper
was motivated by the problem of estimating the covariance matrix of
$n$ independent identically distributed $p$-vectors, $\X_1,\ldots,
\X_N$ with a common $N(\0, \Sigma_p)$ distribution. In
\cite{BickLev08a,BickLev08b} Bickel and Levina show that covariance
matrices which are approximable by banded matrices could be well
estimated, in the operator norm, by banded empirical covariance
matrices, and accordingly their inverses could be approximated by
the inverses of the estimates above. Conversely, inverse covariance
matrices approximable by banded matrices could be well approximated
by data dependent banded matrices and now the covariance matrices
themselves would be approximated by the inverses of the banded
matrices above. The bounds developed above enable us to approximate
both a covariance and its inverse by banded matrices simultaneously
in a very explicit way.

Specifically consider $\Sigma_p$ as the top left $p \times p$ block of a
banded matrix $\Sigma: \ell^2 \to \ell^2$ with a bounded inverse
$\Sigma^{-1}$. Let $\| \cdot \|$ be the operator norm and put
\[
\delta_k (\Sigma) \ :=\ \dist(\Sigma,\BO_k)\ =\ \min \{ \| \Sigma - B \|: B \in \BO_k \}
\ =:\ \| \tilde{B}_k(\Sigma) -\Sigma \|.
\]
with $\tilde{B}_k(\Sigma)\in\BO_k$. Theorem \ref{th:BDO}, specialized to A positive definite
and self-adjoint, says that
\be
\label{eq26}
\delta_{nk} (\Sigma^{-1}) \ \leq\ \frac{2 \delta_k}{m(m-2
\delta_k)} \ +\ \frac{1}{m} \left( \frac{\kappa-1}{\kappa+1} \right)^n,
\ee
where $\delta_k = \delta_k (\Sigma)$, $m^{-1} = m^{-1}(\Sigma)$ is the
norm of $\Sigma^{-1}$, $M=M(\Sigma)$ is the norm of $\Sigma$
and $\kappa=\kappa(\Sigma) = \frac{M}{m}$ is the condition number of $\Sigma$.
Let
\[
\hat \Sigma_p \ :=\ \frac{1}{N} \sum^N_{i=1} (\X_i - \bar \X)
(\X_i - \bar \X)^T
\]
be the empirical covariance of $\X_{p\times 1}$.
The individual elements of $\hat \Sigma_p$, $\hat \sigma_{ij}$ approach the
corresponding $\sigma_{ij}$ as $N \to \infty$ with high probability but
$\hat \Sigma_p $ fails to have an inverse if $p \geq N$ and its
eigenstructure is, in general, diverging from that of $\Sigma$, if $p$
is commensurate or much larger than $N$ i.e. as $p \to \infty$ as well as
$N$ with $\frac{p}{N} \to c$, $0 <c \leq \infty$, see Johnstone \cite{John}.
However, we can, under very mild conditions on $p$ and $N$,
find $k_N \to \infty$ such that
\be
\label{eq27}
\| B_{k_N} (\hat \Sigma_p) - B_{k_N} (\Sigma_p) \| \ \to\ 0
\ee
in probability, where
$$
B_k(A) \ :=\ ( a_{ij}\, \chi(|i-j| \leq k))_{i,j=1}^p
$$
if $A = (a_{ij})_{i,j=1}^p$. In particular, for Gaussian $\X$ this
is true if $\frac{\log p}{N} \to 0$. Therefore if $\Sigma \in \BDO$,
so that $\delta_k(\Sigma) \to 0$ as $k \to \infty$, this yields an
operator norm consistent estimate of $\Sigma$, that is \be
\label{eq28} \| B_{k_N}(\hat \Sigma_p) - \Sigma_p \| \ \leq\ \|
B_{k_N} (\hat \Sigma_p) - B_{k_N} (\Sigma_p) \| \ +\ \| B_{k_N}
(\Sigma_p) - \Sigma_p \| \ \mathop{\longrightarrow}^P\  0. \ee We
can deduce from \eqref{eq26} that if $\| \Sigma^{-1} \| < \infty$
then \be \label{eq29} \| B_{n\,k_N} (\Sigma^{-1}_p) - \Sigma^{-1}_p
\| \ \to\ 0 \ee as $n \to \infty$ slowly with $k_N$. It is, of
course, clear from \eqref{eq27} and \eqref{eq29} that $[  B_{k_N}
(\hat \Sigma_p) ]^{-1}$ will eventually exist, be self-adjoint
positive definite and consistently estimate $\Sigma^{-1}$. However,
this is rather unsatisfactory in practice as well as theory since $
B_{k_N}^{-1}(\hat \Sigma_p)$ in general does not have a band
structure (meaning that it is supported on all its diagonals) in
particular if we assume that $\Sigma^{-1} \in \BO$, so that by
Theorem \ref{th:BDO}, $\Sigma \in \BDO$, our estimate would not
reflect this information. The assumption that $\Sigma^{-1}_p$
belongs to $\BO_k$ has, in the Gaussian case, a statistical
interpretation. It implies that $\X_i$ is independent of $\{ \X_j:
|i-j| >k \}$ given $\{ \X_j: |i-j| \leq k, \ j\neq i \}$. The
assumption that $\Sigma_p \in \BO_k$ has a different interpretation,
implying that $\X_i$ is independent (unconditionally) of $\{ \X_j:
|i-j| > k \}$. One of the interesting consequences of Theorem
\ref{th:BDO} is that it tells us that conditional independence (for
a band structure) cannot occur unless there is approximate
conditional independence and vice versa. That point aside, we are
left with a good but not ``natural'' estimate for $\Sigma^{-1}_p$ if
we assume $\Sigma^{-1} \in \BO$. However, our approach to Theorem
\ref{th:BDO} tells us precisely what to do.

Suppose $B_k(\hat \Sigma_p)$ is our estimate
of $\Sigma_p$. Let $\hat m_k$, $\hat M_k$ be the minimal and maximal absolute
eigenvalues of $B_k(\hat \Sigma_p)$ and
$$
\hat \gamma \ :=\ 2(\hat M_k + \hat m_k)^{-1}.
$$
Then an $nk$ banded estimate of $\Sigma^{-1}_p$ is just
$$
\hat \gamma\ \sum^n_{j=0} \big(I - \hat \gamma B_k(\hat \Sigma_p)
\big)^j.
$$
We are left with the problem of how to choose $k$ and $n$. In theory, if we
have some notion or make assumptions about the magnitude of the ``bias''
$\delta_k(\Sigma)$ and calculate stochastic bounds on the ``variance''
$\| B_k (\hat \Sigma_p) -  B_k (\Sigma_p) \|$ and make assumptions
about how many zeros $\Sigma^{-1}$ has, we can use \eqref{eq26} to estimate
the optimal choices of $k$ and $n$. In practice, it is better to use some
data determined choice, e.g. by crossvalidation, see Bickel, Levina
\cite{BickLev08a,BickLev08b}.
However, we believe that \eqref{eq26} and \eqref{eq28} can be used to compute
minmax bounds and oracle inequalities on the performance of estimates
of $\Sigma^{-1}_p$, from the ones obtained for estimates of $\Sigma_p$, see
Cai, Zhang, Zhou \cite{CaiZhangZhou}. It should be clear that whatever we
have said of banding applies to generalized banding up to a permutation also.

The application of extended generalized banding is satisfactory if
we think of the coordinates of $\X$ as corresponding to labelled points on
a manifold, as is reasonable in (say) geophysical applications, where $\X$ is
the state of some variable such as pressure, at a grid of points on the globe
at some time. But it is not relevant if the labels are meaningless as in
microarrays genomics, where the coordinates simply label genes. However, it
should now be clear that our second generalization to generalized banding up
to some unknown permutation $\pi$ deals with such situations. It enables us
to define classes of matrices $\Sigma$ in terms of their dependency graph defined
in terms of $\Sigma^{-1}$ by having vertices correspond to coordinates of
$\X$ with an edge between $i$ and $j$ if the $(i,j)$th entry of $\Sigma^{-1}$
is different from 0. However, although the determination of $\rh$ is
dictated by the situation, estimation of $\pi$ is nontrivial and will
be pursued elsewhere.

\bigskip
{\bf Acknowledgements. } This research was partially supported by
NSF grants DMS-0605236 and DMS-0906808 to Bickel and by Marie Curie
Grants MEIF-CT-2005-009758 and PERG02-GA-2007-224761 to Lindner.
Both authors have been supported by the Statistics and Applied
Mathematics Institute (NSF) and the American Institute of
Mathematics (AIM).


\bigskip

\noindent {\bf Author's addresses:}\\[0mm]

\noindent
Peter~J.~Bickel\hfill{\tt bickel@stat.berkeley.edu}\\
Department of Statistics\\
367 Evans Hall\\
Berkeley, CA, 94710-3860\\
USA\\[2mm]

\noindent
Marko Lindner\hfill {\tt marko.lindner@mathematik.tu-chemnitz.de}\\
Fakult\"at Mathematik\hfill  (corresponding author)\\
TU Chemnitz\\
D-09107 Chemnitz\\
GERMANY

\end{document}